\documentclass[12]{amsart}
\usepackage{latexsym}
\usepackage{amssymb}
\usepackage{amsfonts, amsmath}
\usepackage[T1]{fontenc}
\usepackage[utf8]{inputenc}
\usepackage{newunicodechar}
\newunicodechar{−}{-}
\usepackage{amsmath, amssymb, amsthm, latexsym}

\newtheorem{lm}{Lemma}
\newtheorem{Thm}{Theorem}

\newtheorem*{corr}{Corollary}

\begin{document}
\title{Unbranched Riemann domains over Stein spaces and Cartier divisors}
\author{}
\noindent
\maketitle
\begin{center}
\em{Youssef Alaoui}\\
\em{y.alaoui@iav.ac.ma}\\
\end{center}
\noindent
\em{Department of Mathematics,
Hassan II Institute of Agronomy}\\ 
\em{and Veterinary Sciences,
Madinat Al Irfane, BP 6202, Rabat, 10101, Morocco,}\\
\date{}
\linespread{1.3}
\date{}

\maketitle
\begin{abstract}
It is proved that an unbranched Riemann domain $\Pi : X\rightarrow Y$ over a complex Stein space $Y$ of dimension $n\geq 2$
is Stein if and only if $X$ is cohomologically $2$-complete with respect to the structure sheaf ${\mathcal{O}}_{X}$
and every topologically trivial holomorphic line bundle over $X$ is associated to a Cartier divisor.
\end{abstract}
\section{Introduction}
It was shown by Abe \cite{ref1} that if $X$ is an open subset of a Stein manifold $Y$ of dimension $n$
such that $H^{k}(X,{\mathcal{O}_{X}})=0$ for $k\geq 2$, then $X$ is Stein if and only if every topologically trivial
holomorphic line bundle $L$ on $X$ is associated to a Cartier divisor on $X$.\\
\hspace*{.1in}This result has been generalized to unbranched Riemann domains  over Stein manifolds 
by Breaz and Vajaitu in \cite{ref3}.\\
\hspace*{.1in}We recall that a complex space $X$ is called a Riemann domain over $Y$ if there exists an analytic map $\Pi : X\rightarrow Y$ which is locally biholomorphic.\\
\hspace*{.1in}It was also shown in \cite{ref1} that if $(Y,{\mathcal{O}}_{Y})$ is a (not necessarily reduced) Cohen-Macauly Stein space
of pure dimension $n$ and $X$ an open set in $Y$ such that :\\
(i) $H^{p}(X,{\mathcal{O}}_{X})=0$ for all $p\geq 2$,\\
(ii)and for every holomorphic line bundle $L$ on $X$ which is an element of the image of the composition of $\Phi$ of the homomorphisms
$$H^{1}(X, {\mathcal{O}}_{Y}|{X})\stackrel{red^{*}}\rightarrow H^{1}(X, {\mathcal{O}}_{X})\stackrel{e^{*}}\rightarrow H^{1}(X,{\mathcal{O}^{*}_{X}})$$

$L$ is associated to some Cartier divisor on $X$,\\
then $X$ is locally Stein at every point $x\in \partial{X}\setminus Sing(Y)$\\
\\
\hspace*{.1in}The main purpose of this article is to generalize the above results to arbitrary complex spaces. More precisely we prove :
\begin{Thm}
Let $\Pi: X \rightarrow Y$ be an unbranched Riemann domain over a (not necessarily reduced) Stein complex space $Y$ of dimension $n\geq 2$ such that $X$ is cohomologically 
$2$-complete with respect to the structure sheaf ${\mathcal{O}}_{X}$.  Then $X$ is Stein if and only if any topologically trivial holomorphic line bundle on $X$ is defined by a Cartier divisor on $X$. 
\end{Thm}
We recall that a complex space $X$ is called cohomologically $q$-complete with respect to a coherent analytic sheaf ${\mathcal{F}}$
on $X$ if the cohomology group $H^{p}(X,{\mathcal{F}})$ vanishes for all $p\geq q$.\\
\hspace*{.1in}The proof of the above theorem is inspired in part by techniques used in the author's preprint on the arxiv \cite{ref2}.\\
\hspace*{.1in}In particular, we obtain the interesting result :
\begin{corr}
{Let $X$ be a Stein space of dimension $n\geq 2$ and $\Omega\subset X$ an
open subset of $X$ such that $H^{p}(\Omega,{\mathcal{O}}_{\Omega})=0$ for all $p\geq 2$. Then $\Omega$ is Stein if and only if every topologically
trivial line bundle on $\Omega$ can be defined by a Cartier divisor on $\Omega$.}
\end{corr}
\section{Preliminaries}
Let $\Omega$ be an open set in $\mathbb{C}^{n}$
with complex coordinates $z_{1},\cdots, z_{n}$. Then 
a function $\phi\in C^{\infty}(\Omega)$ is $q$-convex if its 
Levi form
$$L_{z}(\phi,\xi)=\displaystyle\sum_{i,j}\frac{\partial^{2}\phi(z)}{\partial{z_{i}}\partial{\overline{z}_{j}}}\xi_{i}\overline{\xi}_{j}, \ \ \xi\in \mathbb{C}^{n}$$
has at most $q-1$ negative or zero eingenvalues at any point $z\in\Omega$.\\
\hspace*{.1in}A smooth real valued function $\phi$ on a complex space $X$
is called $q$-convex if every point $x\in X$ has a local chart
$U\rightarrow D\subset \mathbb{C}^{n}$ such that $\phi|_{U}$ has an extension $\hat{\phi}\in C^{\infty}(D, \mathbb{R})$ which is $q$-convex on $D$.\\
\hspace*{.1in}A complex space $X$ is called $q$-complete if there exists a smooth function $\phi : X\rightarrow \mathbb{R}$ which is $q$-convex and
exhaustive on $X$ i.e. $\{\phi<c\}\subset\subset X$ for any $c\in\mathbb{R}$.\\
\hspace*{.1in}The space $X$ is said to be cohomologically $q$-complete
if for every coherent analytic sheaf ${\mathcal{F}}$ on $X$ the cohomology groups $H^{p}(X, {\mathcal{F}})$ vanish for all $p\geq q$.\\
\hspace*{.1in}It is well knnown from the theory of Andreotti and Grauert $[3]$ that if $X$ is $q$-complete, then $X$ is cohomologically $q$-complete. 
\section{Domains over reduced Stein spaces}
\begin{lm}Let $\Pi: X \rightarrow Y$ be an unbranched Riemann domain over a reduced Stein space $Y$ 
of dimension $2$ such that any topologically trivial holomorphic line bundle on $X$ is defined by a Cartier divisor. 
Then $X$ also is Stein.
\end{lm}
In order to prove lemma $1$ we shall need the following theorem due to Breaz and Vajaitu \cite{ref3}.
\begin{Thm}
Let $\Pi : X\rightarrow Y$ be an unbranched Riemann domain over a Stein manifold Y of pure dimension $n\geq 2$ such that $H^j(\mathrm{X}, \mathcal{O})=0$ for $2 \leq j<n$.
Then X is Stein if and only if any topologically trivial holomorpic line bundle on X is defined by a Cartier divisor on X .
\end{Thm}
\begin{proof} We may assume that $Y$ is connected. Then there exists a holomorphic function $f\in {\mathcal{O}(Y)}$
such that its zero set $Z=\{z\in Y : f(z)=0\}$ is a closed analytic subset of pure dimension $1$ and $Sing(Y)\subset Z$.
Moreover, the holomorphic map $\Pi|_{Z'} :Z'=\Pi^{-1}(Z)\rightarrow Z$ is an unbranched Riemann domain over the Stein analytic set $Z$, it follows 
that $Z'$ itself is a Stein , since $Z'$ is of dimension $1$ and every compact analytic subset of $Z'$ is finite.\\
\hspace*{.1in}Since $Z^{\prime}$ is a closed Stein analytic subset of $X$, then, by a theorem of Siu \cite{ref6}, $Z^{\prime}$ admits a Stein open neighborhood $\Omega$ in $X$. Because $\Omega \backslash Z^{\prime}$ also is Stein, then from the exact sequence of local cohomology :
$$
\cdots \rightarrow H^1\left(\Omega \backslash Z^{\prime}, \mathcal{O}_{\Omega}\right) \rightarrow H_{Z^{\prime}}^{2}\left(\Omega, \mathcal{O}_{\Omega}\right) \rightarrow H^{2}\left(\Omega, \mathcal{O}_{\Omega}\right) \rightarrow \cdots
$$
it follows that $H_{Z^{\prime}}^2\left(X, \mathcal{O}_X\right) \cong H_{Z^{\prime}}^2\left(\Omega, \mathcal{O}_{\Omega}\right)=0$ and using the exact sequence
$$
\cdots H^1\left(X, \mathcal{O}_X\right) \rightarrow H^1\left(X \backslash Z^{\prime}, \mathcal{O}_X\right) \rightarrow H_{Z^{\prime}}^{2}\left(X, \mathcal{O}_X\right) \rightarrow \cdots
$$
we see that the restriction
$$
H^1\left(X, \mathcal{O}_X\right) \xrightarrow{r_1} H^1\left(X \backslash Z^{\prime}, \mathcal{O}_X\right)
$$
is surjective.

Consider now the commutative diagrams with exact rows :

$$
\begin{array}{clclcl}
H^1\left(X, \mathcal{O}_X\right) & \xrightarrow{\iota_X} & H^{1}(X, {\mathcal{O}}_{X}^{*}) & \xrightarrow{c_{1, X}} & H^2(X, \mathbb{Z})  \\
\downarrow r_1 & & \downarrow r_2 & & \downarrow & \\
H^1\left(X \backslash Z^{\prime}, \mathcal{O}_{X \backslash Z^{\prime}}\right) & \xrightarrow{\iota_{X \backslash Z^{\prime}}} & \operatorname H^{1}(X\setminus Z^{\prime}, {\mathcal{O}}^{*}_{X}) & \xrightarrow{c_{1, X \backslash Z^{\prime}}} & H^2\left(X \backslash Z^{\prime}, \mathbb{Z}\right) & \rightarrow 0
\end{array}
$$
Let $L$ be a line bundle on $X \backslash Z^{\prime}$ which is topologically trivial. Then there exists a cohomology class $\xi \in H^1\left(X \backslash Z^{\prime}, \mathcal{O}_{X \backslash Z^{\prime}}\right)$ such that $\iota_{X \backslash Z^{\prime}}(\xi)=(L) \in H^1\left(X \backslash Z^{\prime}, \mathcal{O}^{*}_{X \backslash Z^{\prime}}\right)$, where $(L)$ is the equivalence class of cohomology in $H^{1}(X\setminus Z^{\prime}, {\mathcal{O}}_{X}^{*})\cong Pic(X\setminus Z^{\prime})$ defined by $L$. Since the restriction map $r_1: H^1\left(X, \mathcal{O}_X\right) \rightarrow H^1\left(X \backslash Z^{\prime}, \mathcal{O}_{X \backslash Z^{\prime}}\right)$ is surjective, there exists a cohomology class $\xi_1 \in H^1\left(X, \mathcal{O}_X\right)$ such that $r_1\left(\xi_1\right)=\xi$.\\
Therefore, $(\tilde{L}) = \iota_X(\xi_1) \in H^{1}(X,\mathcal{O}_X^{*})$ 
is the equivalence class of a line bundle that is topologically trivial on $X$, 
and whose restriction to $X \setminus Z$ coincides with $L$.
By assuption, $\tilde{L}=O_X(D)$ for some Cartier divisor $D$ on $X$. Then $\left.\tilde{L}\right|_{X \backslash Z^{\prime}}=L=O_{X \backslash Z^{\prime}}\left(\left.D\right|_{X \backslash Z^{\prime}}\right)$.
Since $\Pi|_{X'} : X'=X\setminus Z'\rightarrow Y\setminus Z$ is an unbranched Riemann domain over the Stein manifold $Y\setminus Z$ such that
every topologically trivial holomorphic line bundle on  $X\setminus Z'$ is defined by a Cartier divisor, and the cohomology group 
$H^{p}(X \backslash Z^{\prime}, {\mathcal{O}}_{X})$ vanishes for all $p\geq 2$, $X'$ being of dimension $2$ and obviously 
has no compact analytic subset of positive dimension, it follows from Theorem $1$ that $X'$ is a Stein manifold.\\

\hspace*{.1in}Let now $\xi : \tilde{X}\rightarrow X$ be a resolution of singularities,
i.e. $\tilde{X}$ is a complex manifold and $\xi$ is a proper modification such that the induced
map
$$\tilde{X}\backslash\{\xi^{-1}(Sing(X))\}\rightarrow X\setminus Sing(X)\}$$
is biholomorphic. Let ${\mathcal{F}}$ be a coherent analytic sheaf on $X$. Then there exists a canonical sheaf homomorphism
${\mathcal{F}}\stackrel{\psi}\rightarrow \xi_{*}\xi^{*}({\mathcal{F}}).$
If we set ${\mathcal{F}}_{1}=ker \ \psi$ and ${\mathcal{F}}_{2}=Im \ \psi$,
then $Supp ({\mathcal{F}}_{1})\subset Z'$ and there is an exact sequence
$$0\rightarrow {\mathcal{F}}_{1}\rightarrow {\mathcal{F}}\stackrel{\psi}\rightarrow {\mathcal{F}}_{2}\rightarrow 0$$
Let ${\mathcal{I}}(Z')\subset {\mathcal{O}}_{X}$ be the subsheaf of germs of holomorphic functions which vanish on $Z'$. Let ${\mathcal{O}}_{Z'}={\mathcal{O}}_{X}/_{{\mathcal{I}}(Z')}$
and $({\mathcal{F}}_{2})_{Z'}={\mathcal{F}}_{2}\otimes_{{\mathcal{O}}_{X}}{\mathcal{O}}_{Z'}$.
If $e$ is the image in ${\mathcal{O}}_{Z'}$ of the section $1$ on ${\mathcal{O}}_{X}$,
then any element of $(({\mathcal{F}}_{2})_{Z'})_{x}$ can be written in the form
$\xi\otimes e_{x},$
where $\xi\in ({\mathcal{F}}_{2})_{x}$. Then the homomorphism
$\eta : {\mathcal{F}}_{2}\rightarrow ({\mathcal{F}}_{2})_{Z'}$ defined by
$\eta(\alpha)=\alpha\otimes e$ is surjective and we have the exact sequence
$$0\rightarrow Ker (\psi)\rightarrow Ker (\eta o\psi)\rightarrow \frac{Ker (\eta o\psi)}{Ker (\psi)}\rightarrow 0$$
Since clearly $Supp \ Ker \ \psi\subset Sing(X)\subset Z'$ and
$Supp \ \frac{Ker (\eta o\psi)}{Ker (\psi)}\subset X\setminus Z'$, 
then $H^{1}(Z', Ker \ \psi)=H^{1}(X\setminus Z', \frac{Ker \ \eta o\psi}{Ker \ \psi})=0$,
Therefore, from the long exact sequence of cohomology
$$\cdots\rightarrow H^{1}(Z', Ker \ \psi)\rightarrow H^{1}(X, Ker \ \eta o\psi )
\rightarrow H^{1}(X\setminus Z', \frac{Ker \ \eta o\psi}{Ker \ \psi})\rightarrow\cdots$$
we deduce that $H^{1}(X, Ker \ \eta o\psi )=0$.\\
Moreover, since $H^{1}(X,({\mathcal{F}}_{2})_{Z'})=H^{1}(Z',({\mathcal{F}}_{2})_{Z'})=0$, then
by using the long exact sequence of cohomology associated to the exact sequence of sheaves
$$0\longrightarrow Ker (\eta o\psi)\longrightarrow {\mathcal{F}}\stackrel{\eta o\psi}\longrightarrow
({\mathcal{F}}_{2})_{Z'}\longrightarrow 0$$
one obtains $H^{1}(X, {\mathcal{F}})=0$, which completes the proof of lemma $1$.
\end{proof}
\begin{lm}Let $\Pi: X \rightarrow Y$ be an unbranched Riemann domain over a reduced Stein space $Y$ 
of dimension $n\geq 2$ such that the cohomology group $H^{p}(X, {\mathcal{O}_{X}})=0$ for $p\geq 2$
and any topologically trivial holomorphic line bundle on $X$ is defined by a Cartier divisor. Then $X$ also is Stein.
\end{lm}
\begin{proof}
That $X$ is Stein was proved in the case $n=2$ in the proof of lemma $1$. We now suppose that $n\geq 3$
and that lemma $2$ has already been proved in dimension $\leq n-1$.
We may, of course assume that $Y$ is connected and let $f$ be a holomorphic function $f\in {\mathcal{O}}(Y)$ such that $Z=\{y\in Y : f(y)=0\}$ 
is of pure dimension $n-1$ and $Sing(Y)\subset Z$.\\
\hspace*{.1in}The multiplication by $f o\Pi$ defines an injective morphism :
$$0\longrightarrow {\mathcal{O}_{X}}\longrightarrow {\mathcal{O}_{X}}$$
and, if $Z'=\Pi^{-1}(Z)$, then from the exact sequence of sheaves :
$$0\longrightarrow {\mathcal{O}_{X}}\longrightarrow {\mathcal{O}_{X}}\longrightarrow {\mathcal{O}_{X}}/_{(f o\Pi){\mathcal{O}_{X}}}\rightarrow 0$$
we deduce the long exact sequence of cohomology :
$$\cdots\rightarrow H^{p}(X, {\mathcal{O}_{X}})\rightarrow H^{p}(Z', {\mathcal{O}_{Z'}})\rightarrow 
H^{p+1}(X, {\mathcal{O}_{X}})\rightarrow\cdots$$
Since $H^{p}(X, {\mathcal{O}_{X}})=0$ for all $p\geq 2$, it follows that $H^{p}(Z', {\mathcal{O}_{Z'}})=0$ for $p\geq 2$
and the restriction map 
$$H^{1}(X, {\mathcal{O}_{X}})\stackrel{r_{1}}\rightarrow H^{1}(Z', {\mathcal{O}_{Z'}})$$
is surjective.\\ 
\hspace*{.1in}Now, from the exponential exact sequence of sheaves :
$$0 \rightarrow \mathbb{Z}_{Z'} \xrightarrow{\times2 \pi i} \mathcal{O}_{Z'} \xrightarrow{\exp ()} \mathcal{O}_{Z'}^{*} \rightarrow 0$$
One obtains the commutative diagrams of restriction maps with exact rows :
$$
\begin{aligned}
&\begin{array}{cccccc}
H^1(X,\mathcal{O}_X) & \xrightarrow{\ \iota_X\ } & \operatorname{Pic}(X) & \xrightarrow{\ c_{1,X}\ } & H^2(X,\mathbb{Z}) & \to 0 \\[0.4em]
\downarrow r_{1} & & \downarrow r_{2} & & \downarrow & \\[0.4em]
H^1(Z',\mathcal{O}_{Z'}) & \xrightarrow{\ \iota_Z'\ } & \operatorname{Pic}(Z') & \xrightarrow{\ c_{1,Z'}\ } & H^2(Z',\mathbb{Z}) & \to 0
\end{array}
\end{aligned}
$$
Let $L$ be a line bundle on $Z'$ which is topologically trivial. Then there exists a cohomology class $\xi\in H^1(Z',\mathcal{O}_{Z'})$
such that $\iota_{Z'}(\xi)=(L)\in Pic(Z')$. Since the restriction map $r_{1} : H^1(X,\mathcal{O}_X)\rightarrow H^1(Z',\mathcal{O}_{Z'})$ is surjective,
there exists a cohomology class $\xi_{1}\in H^1(X,\mathcal{O}_X)$ such that $r_{1}(\xi_{1})=\xi$. Therefore,
$(\tilde{L})=\iota_{X}(\xi_{1})\in Pic(X)$ is an equivalent class of line bundles which are  topologically trivial on $X$ and whose restriction to $Z'$ coincide with $L$. By assuption, $\tilde{L}=O_{X}(D)$ for some Cartier divisor $D$ on $X$. Then $\tilde{L}|_{Z'}=L=O_{Z'}(D|_{Z'})$. As
the restriction map $\Pi|_{Z'} : Z'\rightarrow Z$ is an unbranched Riemann domain over the Stein space $Z$ such that $H^{p}(Z', {\mathcal{O}}_{Z'})=0$
for all $p\geq 2$ and any topologically holomorphic line bundle $L$ on $Z'$ is associated to a Cartier divisor on $Z'$, then by the induction 
hypothesis it follows that $Z'$ is a Stein space.\\
\hspace*{.1in}Since $Z'$ is a closed Stein analytic subset of $X$, then $Z'$ admits a Stein open neighborhood $\Omega$
in $X$. Therefore by using the exact sequence of cohomology
$$\cdots\rightarrow H^{p}(\Omega\setminus Z', {\mathcal{O}}_{\Omega})\rightarrow H^{p+1}_{Z'}(\Omega, {\mathcal{O}}_{\Omega})\rightarrow
H^{p+1}(\Omega, {\mathcal{O}}_{\Omega})\rightarrow\cdots$$
and noting that $\Omega$ and $\Omega\setminus Z'$ are Stein, we find that
$H^{p}_{Z'}(X, {\mathcal{O}}_{X})\cong H^{p}_{Z'}(\Omega, {\mathcal{O}}_{\Omega})=0$ for all $p\geq 2$ and from the exact sequence of local cohomology
$$\cdots H^{p}(X, {\mathcal{O}}_{X})\rightarrow H^{p}(X\setminus Z', {\mathcal{O}}_{X})\rightarrow H^{p+1}_{Z'}(X, {\mathcal{O}}_{X})\rightarrow\cdots$$
it follows that $H^{p}(X\setminus Z', O_{X})=0$ for all $p\geq 2$ and 
the restriction 
$$H^{1}(X, {\mathcal{O}}_{X})\stackrel{r_{1}}\rightarrow H^{1}(X\setminus Z', {\mathcal{O}}_{X})$$
is surjective.\\
\hspace*{.1in}Consider now the commutative diagrams with exact rows :
\[
\begin{aligned}
&\begin{array}{cccccc}
H^1(X,\mathcal{O}_X) 
& \xrightarrow{\ \iota_X\ } 
& \operatorname{Pic}(X) 
& \xrightarrow{\ c_{1,X}\ } 
& H^2(X,\mathbb{Z}) 
& \to 0 \\[0.4em]
\downarrow r_{1} 
& 
& \downarrow r_{2} 
& 
& \downarrow 
& \\[0.4em]
H^1(X\setminus Z',\mathcal{O}_{X\setminus Z'}) 
& \xrightarrow{\ \iota_{X\setminus Z'}\ } 
& \operatorname{Pic}(X\setminus Z') 
& \xrightarrow{\ c_{1,X\setminus Z'}\ } 
& H^2(X\setminus Z',\mathbb{Z}) 
& \to 0
\end{array}
\end{aligned}
\]
Let $L$ be a line bundle on $X\setminus Z'$ which is topologically trivial. Then a similar proof as that of lemma $1$ shows that $L$
is associated to a Cartier divisor on $X\setminus Z'$. Since in addition
the restriction map $\Pi|_{X\setminus Z'} : X\setminus Z'\rightarrow Y\setminus Z$ is an unbranched Riemann domain over the Stein manifold 
$Y\setminus Z$ such that $H^{p}(X\setminus Z', {\mathcal{O}}_{X\setminus Z'})=0$
for all $p\geq 2$, then by Theorem $1$ it follows that $X\setminus Z'$ is a Stein manifold.\\
\hspace*{.1in}Let now ${\mathcal{F}}$ be a coherent analytic sheaf on $X$, and let $\xi: \tilde{X} \rightarrow X$ be a resolution of singularities.
Then there exists a canonical sheaf homomorphism
\[
\mathcal{F} \xrightarrow{\psi} \xi_* \xi^*(\mathcal{F}).
\]
If we set $\mathcal{F}_1 = \operatorname{Ker} \psi$ and $\mathcal{F}_2 = \operatorname{Im} \psi$, then $\operatorname{Supp}(\mathcal{F}_1)\subset Z'$, and we obtain the short exact sequence
\[
0 \longrightarrow \mathcal{F}_1 \longrightarrow \mathcal{F} \xrightarrow{\psi} \mathcal{F}_2 \longrightarrow 0.
\]
Let $\mathcal{I}(Z') \subset \mathcal{O}_X$ denote the subsheaf of germs of holomorphic functions vanishing on $Z'$. Define $\mathcal{O}_{Z'} = \mathcal{O}_X / \mathcal{I}(Z')$ and $(\mathcal{F}_2)_{Z'} = \mathcal{F}_2 \otimes_{\mathcal{O}_X} \mathcal{O}_{Z'}$.
If $e$ is the image in $\mathcal{O}_{Z'}$ of the unit section $1 \in \mathcal{O}_X$, then any element of $\left((\mathcal{F}_2)_{Z'}\right)_x$ can be expressed as $\xi \otimes e_x$, where $\xi \in (\mathcal{F}_2)_x$.

Define a homomorphism
\[
\eta: \mathcal{F}_2 \longrightarrow (\mathcal{F}_2)_{Z'}
\]
by $\eta(\alpha) = \alpha \otimes e$.
This map is surjective, and we have the exact sequence
\[
0 \longrightarrow \operatorname{Ker}(\psi) \longrightarrow \operatorname{Ker}(\eta \circ \psi) \longrightarrow
\frac{\operatorname{Ker}(\eta \circ \psi)}{\operatorname{Ker}(\psi)} \longrightarrow 0.
\]

Since it is clear that $\operatorname{Supp}(\operatorname{Ker} \psi) \subset \operatorname{Sing}(X) \subset Z'$ and
\[
\operatorname{Supp}\left(\frac{\operatorname{Ker}(\eta \circ \psi)}{\operatorname{Ker}(\psi)}\right) \subset X \setminus Z',
\]
the associated long exact sequence of cohomology 
\[
\cdots \longrightarrow H^1(Z', \operatorname{Ker} \psi) \longrightarrow
H^1(X, \operatorname{Ker}(\eta \circ \psi)) \longrightarrow
H^1(X \setminus Z', \tfrac{\operatorname{Ker}(\eta \circ \psi)}{\operatorname{Ker} \psi}) \longrightarrow \cdots
\]
implies that $H^1(X, \operatorname{Ker} \eta o \psi)=0$. Furthermore, since $H^1\left(X,\left(\mathcal{F}_2\right)_{Z^{\prime}}\right)= H^1\left(Z^{\prime},\left(\mathcal{F}_2\right)_{Z^{\prime}}\right)=0$, then by using the long exact sequence of cohomology associated to the exact sequence of sheaves
$$
0 \longrightarrow \operatorname{Ker}(\eta o \psi) \longrightarrow \mathcal{F} \xrightarrow{\eta o \psi}\left(\mathcal{F}_2\right)_{Z^{\prime}} \longrightarrow 0
$$
we see that $H^1(X, \mathcal{F})=0$, which shows that $X$ is Stein.
\end{proof}
\section{Proof of theorem $1$}

Assume now that $Y$ is not necessarily reduced. 
A direct verification shows that the reduction map
$\Pi_{red}:X_{red}\rightarrow Y_{red}$ is an unbranched Riemann domain with the same properties as the original morphism 
$\Pi : X\rightarrow Y$ of theorem $1$. 
Indeed, let $x \in X$. By assumption, there exist open neighborhoods
\[
x \in U \subset X, \quad \Pi(x) \in V \subset Y
\]
such that
\[
\Pi|_U : (U, \mathcal O_X) \xrightarrow{\;\sim\;} (V, \mathcal O_Y)
\]
is an isomorphism of complex spaces.  

Taking reductions preserves isomorphisms of ringed spaces, so
\[
(U, \mathcal O_{U_{\mathrm{red}}}) \xrightarrow{\;\sim\;} (V, \mathcal O_{V_{\mathrm{red}}})
\]
is also an isomorphism.  

By definition, this shows that
\[
\Pi_{\mathrm{red}} : X_{\mathrm{red}} \longrightarrow Y_{\mathrm{red}}
\]
is locally biholomorphic, and hence an unbranched Riemann domain.\\
We are now going to show that
\[
H^p(X, \mathcal{O}_X) \cong H^p(X, \mathcal{O}_{X_{\mathrm{red}}}) \quad \text{for all } p \ge 0.
\]
Let $\mathcal{N} \subset \mathcal{O}_X$ denote the nilradical sheaf, so that
\[
0 \longrightarrow \mathcal{N} \longrightarrow \mathcal{O}_X \longrightarrow \mathcal{O}_{X_{\mathrm{red}}} \longrightarrow 0
\]
is exact. Since $\mathcal{N}$ is nilpotent, there exists $k \ge 1$ such that $\mathcal{N}^k = 0$. By induction on the powers of $\mathcal{N}$, one shows that
\[
H^p(X, \mathcal{N}) = 0 \quad \text{for all } p \ge 1.
\]
The long exact sequence in cohomology then gives isomorphisms
\[
H^p(X, \mathcal{O}_X) \xrightarrow{\sim} H^p(X, \mathcal{O}_{X_{\mathrm{red}}}) \quad \text{for all } p \ge 0.
\]

It remains to prove that if every topologically trivial holomorphic line bundle on $X$
is associated to a Cartier divisor, then the same property holds on $X_{\mathrm{red}}$.
In fact, let $L_{\mathrm{red}}$ be a topologically trivial holomorphic line bundle on $X_{\mathrm{red}}$.
Since the sheaf of meromorphic functions is insensitive to nilpotents, the
reduction morphism induces a canonical isomorphism of Cartier divisor groups
\[
\mathrm{Div}(X_{\mathrm{red}})\;\xrightarrow{\;\sim\;}\;\mathrm{Div}(X).
\]
Moreover, it induces an isomorphism
\[
\operatorname{Pic}(X_{\mathrm{red}})\;\xrightarrow{\;\sim\;}\;\operatorname{Pic}(X).
\]
Consider the commutative diagram
\[
\begin{array}{ccc}
\mathrm{Div}(X_{\mathrm{red}}) & \xrightarrow{\;\;\mathcal O_{X_{red}}\;\;} &
\operatorname{Pic}(X_{\mathrm{red}}) \\
\Big\downarrow{\alpha}_{\simeq} & & \Big\downarrow{\beta}_{\simeq} \\
\mathrm{Div}(X) & \xrightarrow{\;\;\mathcal O_{X}\;\;} &
\operatorname{Pic}(X),
\end{array}
\]
where $\alpha$ is the canonical identification of Cartier divisors induced by
the isomorphism of meromorphic function sheaves, and $\beta$ is the canonical
isomorphism of Picard groups induced by the reduction morphism.

Let $L_{\mathrm{red}} \in \operatorname{Pic}(X_{\mathrm{red}})$ be topologically
trivial. By hypothesis, there exists a Cartier divisor $D \in \mathrm{Div}(X)$
such that
\[
\mathcal O_X(D) \simeq \beta(L_{\mathrm{red}}).
\]
Set $D_{\mathrm{red}} := \alpha^{-1}(D) \in \mathrm{Div}(X_{\mathrm{red}})$.

By commutativity of the diagram, we have
\[
\beta\big(\mathcal O_{X_{\mathrm{red}}}(D_{\mathrm{red}})\big)
=
\mathcal O_X(D)
\simeq
\beta(L_{\mathrm{red}}).
\]
Since $\beta$ is injective, it follows that
\[
L_{\mathrm{red}} \simeq \mathcal O_{X_{\mathrm{red}}}(D_{\mathrm{red}}).
\]
Since $X_{red}$ is Stein by the previous demonstration, a theorem of Cartan and Serre implies that $X$
itself is Stein. Thus theorem $1$ is completely proved.



\end{document}